\documentclass[11pt]{article}
\usepackage{latexsym,amssymb}
\usepackage{psfrag}
\usepackage{amsmath}
\usepackage{color}
\usepackage{slashbox}
\usepackage{lineno}

\newtheorem{Theorem}{\bf Theorem}[section]
\newtheorem{Lemma}{\bf Lemma}[section]
\newtheorem{Proposition}{\bf Proposition}[section]
\newtheorem{Corollary}{\bf Corollary}[section]
\newtheorem{Remark}{\bf Remark}[section]
\newtheorem{Example}{\bf Example}[section]
\newtheorem{Definition}{\bf Definition}[section]

\newenvironment{theorem}{\begin{Theorem}$\!\!\!$}{\end{Theorem}}

\newenvironment{proposition}{\begin{Proposition}$\!\!\!$}{\end{Proposition}}

\newenvironment{remark}{\begin{Remark}$\!\!\!$}{\end{Remark}}

\newenvironment{definition}{\begin{Definition}$\!\!\!$}{\end{Definition}}

\def\Xint#1{\mathchoice
{\XXint\displaystyle\textstyle{#1}}%
{\XXint\textstyle\scriptstyle{#1}}%
{\XXint\scriptstyle\scriptscriptstyle{#1}}%
{\XXint\scriptscriptstyle\scriptscriptstyle{#1}}%
\!\int}
\def\XXint#1#2#3{{\setbox0=\hbox{$#1{#2#3}{\int}$}
\vcenter{\hbox{$#2#3$}}\kern-.5\wd0}}

\def\dashint{\Xint-}

\numberwithin{equation}{section}
\newcommand{\R}{{\bf R}}
\newcommand{\RN}{{\bf R}^N}
\newcommand{\RNpls}{{\bf R}^N_+}
\newcommand{\half}{\frac{1}{2}}

\begin{document}
\title {Sharp estimate of the life span of solutions to \\the heat equation with a nonlinear boundary condition}
\author{
Kotaro Hisa 
}
\date{}
\maketitle
\begin{abstract}
Consider the heat equation with a nonlinear boundary condition
\begin{equation*}
{\rm (P)}\qquad
\left\{
\begin{array}{ll}
\partial_t u=\Delta u,\qquad & x\in{\bf R}^N_+,\,\,\,t>0,\vspace{3pt}\\
\displaystyle{-\frac{\partial u}{\partial x_N} u}=u^p, & x\in\partial{\bf R}^N_+,\,\,\,t>0,\vspace{3pt}\\
u(x,0)=\kappa\psi(x),\qquad  & x\in \overline{{\bf R}^N_+},
\end{array}
\right.
\qquad\qquad
\end{equation*}
where $N\ge 1$, $p>1$, $\kappa>0$
and $\psi$ is a nonnegative measurable function in ${\bf R}^N_+ :=\{y\in\RN:y_N>0 \}$. 
Let us denote by $T(\kappa\psi)$ the life span of solutions to problem~(P).
We investigate the relationship between the singularity of $\psi$ at the origin and $T(\kappa\psi)$ for sufficiently large $\kappa>0$
and
the relationship between the behavior of $\psi$ at the space infinity and $T(\kappa\psi)$ for sufficiently small $\kappa>0$.
Moreover, we obtain sharp estimates of $T(\kappa\psi)$, as $\kappa\to\infty$ or $\kappa\to+0$. 
\end{abstract}
\vspace{15pt}
\noindent Address:

\smallskip
\noindent K.~H.:Graduate School of Mathematical Sciences, The University of Tokyo,\\
\qquad\,\,\, 3-8-1 Komaba, Meguro-ku, Tokyo 153-8914, Japan. \\\noindent 
E-mail: {\tt khisa@ms.u-tokyo.ac.jp}\\
\vspace{15pt}
\vspace{3pt}
\newpage

\section{Introduction}
Consider the heat equation with a nonlinear boundary condition
\begin{equation}
\label{prob}
\left\{
\begin{array}{ll}
\partial_t u=\Delta u,\quad & x\in{\bf R}^N_+,\,\,\,t>0,\vspace{3pt}\\
\displaystyle{-\frac{\partial u}{\partial x_N}=u^p}, & x\in\partial{\bf R}^N_+,\,\,\,t>0,\
\end{array}
\right.
\qquad\qquad
\end{equation}
with the initial condition
\begin{equation}
\label{ini}
u(x,0)=\kappa\psi(x),\qquad   x\in D:=\overline{{\bf R}^N_+},
\end{equation}
where $N\ge1$, $p>1$, $\kappa>0$ 
and $\psi$ is a nonnegative measurable function in ${\bf R}^N_+:=\{y\in\RN:y_N>0 \}$.
Let $T(\kappa\psi)$ denote the maximal existence time of the minimal solution to problem~\eqref{prob} with \eqref{ini}.
We call $T(\kappa\psi)$ the life span of solutions to problem~\eqref{prob} with \eqref{ini} (see Definitions~\ref{Definition:1.1} and \ref{Definition:1.2}).
The life spans  depend on a lot of factors such as diffusion effect, nonlinearity of equations, boundary conditions and the singularity or the decay of initial functions
and they have been studied,
see e.g., \cite{FR, HI19, IS1, IS2}.
For  related results on semilinear parabolic equations, see e.g., \cite{D, FR, GW, HI18, HI19, IkSo1, IS1, IS2, LN,  MY1, MY2, OY, S, YY, Y} and references therein.

Problem~\eqref{prob} can be physically interpreted as a nonlinear radiation law 
and it has been studied in many papers (see e.g., \cite{ACR1, AR, DFL, FR, FiloK, GL, HI19, IK, IS1, IS2}).
Among others, the author of this paper and Ishige \cite{HI19} obtained necessary conditions and  sufficient conditions for the solvability of problem~\eqref{prob}
and identified the strongest singularity 
of the initial function for the existence of  solutions to problem~\eqref{prob}.
In this paper,  applying the results in \cite{HI19}, we obtain sharp estimates of the life span $T(\kappa\psi)$ as $\kappa\to\infty$ or $\kappa\to+0$ and show that the behavior of the life span $T(\kappa\psi)$ as $\kappa\to\infty$ and $\kappa\to+0$ depends on the singularity and the decay of $\psi$, respectively. 
The proofs of our results require  careful treatments of parameters in the results in \cite{HI19}.

Before stating the main results of this paper, we have to define the life span $T(\kappa\psi)$ of solutions to \eqref{prob} with \eqref{ini} exactly.
To do that, we formulate  the definition of solutions to \eqref{prob}.
Let $G=G(x,y,t)$ be the Green function for the heat equation on $\RNpls$ with the homogeneous Neumann boundary condition.
For $y=(y_1,\cdots,y_N)\in\RN$, $y'$ is given by $y'=(y_1,\cdots,y_{N-1})$.
\begin{definition}
\label{Definition:1.1}
Let $u$ be a nonnegative and continuous function in $D\times(0,T)$, where $0<T<\infty$. 
\begin{itemize}

\item Let $\varphi$ be a nonnegative measurable function in $\RNpls$. We say that $u$ is a solution to \eqref{prob} in $[0,T)$ with $u(0) = \varphi$ if $u$ satisfies 
\begin{equation*}
\label{IE}
u(x,t) = \int_D G(x,y,t) \varphi(y) \, dy + \int_0^t \int_{{\bf R}^{N-1}} G(x,y',0, t-s) u(y',0,s)^p \, dy' ds
\end{equation*}
for $(x,t) \in D\times (0,T).$

\item  We say that $u$ is a minimal solution to \eqref{prob} in $[0,T)$ with $u(0)=\varphi$ if $u$ is a  solution to \eqref{prob} in $[0,T)$ with $u(0)=\varphi$ and satisfies
 $$
 u(x,t) \le w(x,t) \quad \mbox{in} \quad D\times(0,T)
 $$ 
for any solution $w$ to \eqref{prob} in $[0,T)$ with $w(0)=\varphi$.
\end{itemize}
\end{definition}
\begin{remark}
Let $u$ be a solution to problem~\eqref{prob} with $u(0)=\varphi$ in the sense of Definition~{\rm \ref{Definition:1.1}}. Then $u$ satisfies the initial condition in the sense of distributions, that is,
$$
\lim_{t\to+0} \int_{D}u(y,t)\eta(y) \, dy = \int_{D}\varphi(y)\eta(y)\, dy
$$
for all $\eta\in C_0(\RN)$.
\end{remark}

Since the minimal solution is unique, we can define the life span $T(\kappa\psi)$ as follows:
\begin{definition}
\label{Definition:1.2}
The life span $T(\kappa\psi)$ of  solutions to \eqref{prob} with \eqref{ini} is defined by the maximal existence time of the minimal solution to \eqref{prob} with \eqref{ini}.
\end{definition}

Next, we set up notation.
For any $x\in\RN$ and $r>0$, set 
$$
B_+(x,r) := \{y\in\RN : |x-y|<r\}\cap D.
$$
For any set $E$, let $\chi_E$ be the characteristic function which has value 1 in $E$ and value 0 outside $E$. 
For any two nonnegative functions $f_1$ and $f_2$ defined in $(0,\infty)$, we write $f_1(\tau)\thicksim f_2(\tau)$ as $\tau\to\infty  ({\rm resp.} +0)$ if there exists a constant $C>0$ such that $C^{-1}f_2(\tau) \le f_1(\tau) \le C f_1(\tau)$ for sufficiently large (resp. small) $\tau>0$.

Now we are ready to state the main results of this paper. In Theorem~\ref{Theorem:1.1} we obtain the relationship between the singularity of $\psi$ and the life span $T(\kappa\psi)$ as $\kappa\to\infty$ and give  sharp estimates to the life span as $\kappa\to\infty$.
Appendix contains a brief summary of Theorem~\ref{Theorem:1.1} (see Tables 1, 2 and 3 in Appendix).
In what follows we set $p_*:=1+{1}/{N}$.

\begin{theorem}
\label{Theorem:1.1}
Assume that 
  $$
  \psi(x):=|x|^{-A}\biggr[\log\left(e+\frac{1}{|x|}\right)\biggr]^{-B}\chi_{B_+(0,1)}(x)\in L^1({\bf R}^N_+)
  \setminus L^\infty({\bf R}^N_+), 
  $$
  where $0\le A \le N$ and
\begin{equation}
\label{eq:1.5}
B>0\quad\mbox{if}\quad A=0,
\qquad
B\in{\bf R}\quad\mbox{if}\quad 0<A<N,
\qquad
B>1\quad\mbox{if}\quad A=N.
\end{equation}
  Then $T(\kappa\psi)\to 0$ as $\kappa\to\infty$ and the following holds:
  \begin{itemize}
  \item[\rm{(i)}]
  $T(\kappa\psi)$ satisfies
  $$
  T(\kappa\psi)\thicksim
  \left\{
  \begin{array}{ll}
  \left[\kappa(\log\kappa)^{-B}\right]^{-\frac{2(p-1)}{-A(p-1)+1}}
   & \quad\mbox{if}\quad
   \displaystyle{{A <\min\biggr\{N,\frac{1}{p-1}\biggr\}}},\vspace{5pt}\\
   \left[\kappa(\log\kappa)^{-B+1}\right]^{-\frac{2(p-1)}{{-A}(p-1)+1}}
   & \quad\mbox{if}\quad 1<p<p_*,\,\,{A=N},\,\, B>1,
  \end{array}
  \right.
   $$
   and  
   $$
   |\log T(\kappa\psi)|\thicksim 
   \left\{
  \begin{array}{ll}
  \kappa^{\frac{1}{B}} & \quad\mbox{if}\quad p>p_*,\,\,\displaystyle{A=\frac{1}{p-1}},\,\,B>0,\vspace{5pt}\\
  \kappa^{\frac{1}{B-N-1}} & \quad\mbox{if}\quad p=p_*,\,\,A=N,\,\, B>N+1,
  \end{array}
  \right.
  $$
  as $\kappa\to\infty$;
  \item[\rm{(ii)}]
  Let $p>p_*$. If, either 
  $$
  A>1/(p-1)  \quad \mbox{and} \quad B\in{\bf R} \qquad or \qquad       
  A=1/(p-1) \quad \mbox{and} \quad B<0,
  $$
  then problem~\eqref{prob} with \eqref{ini} possesses no              local-in-time solutions
  for all $\kappa>0$. 
  If 
  $$
  A=1/(p-1)\quad and \quad B=0, 
  $$
  then problem~\eqref{prob} with \eqref{ini} possesses no local-in-time solutions 
  for sufficiently large $\kappa>0$; 
  \item[\rm{(iii)}]
  Let $p=p_*$. If 
  $$
  A=N\quad \mbox{and}\quad B<N+1, 
  $$
  then problem~\eqref{prob} with \eqref{ini} possesses no local-in-time solutions  
  for all $\kappa>0$. 
  If 
  $$
  A=N\quad \mbox{and}\quad B=N+1, 
  $$
  then problem~\eqref{prob} with \eqref{ini} possesses no local-in-time solutions  
  for  sufficiently large $\kappa>0$. 
  \end{itemize} 
\end{theorem}
\vspace{5pt}
We remark that when $\psi$ is as in Theorem~\ref{Theorem:1.1}, $\psi$ satisfies \eqref{eq:1.5} if and only if $\psi\in L^1_{loc}(\RNpls)$.
It is obvious that $T(\kappa\psi)=0$ for all $\kappa>0$ if \eqref{eq:1.5} does not hold.
\begin{remark}
When $B=0$,
Ishige and Sato {\rm\cite{IS1}}  have already obtained sharp estimates of the life span $T(\kappa\psi)$ as $\kappa\to\infty$ in the case when $\psi(x) = |x|^{-A}$
in a neighborhood of the origin, where 
\begin{equation*}
\label{IS}
0\le A < N \quad \mbox{if} \quad 1<p<p_* \qquad \mbox{and} \qquad 0\le  A < \frac{1}{p-1} \quad \mbox{if} \quad p\ge p_*,
\end{equation*}
and proved that $T(\kappa\psi)\thicksim \kappa^{-\frac{2(p-1)}{-A(p-1)+1}}$
as $\kappa\to\infty$. This also follows from  Theorem~{\rm \ref{Theorem:1.1}}.
\end{remark}

\begin{remark}
Let $N=1$ and let $\psi$ be a continuous, positive and bounded function in $\R$.
Fern\'andez Bonder and Rossi {\rm\cite{FR}} obtained the precise asymptotic behavior of the life span $T(\kappa\psi)$ as $\kappa\to\infty$, that is, $\lim_{\kappa\to\infty}\kappa^{2(p-1)} T(\kappa\psi) = T(\psi(0)).$
\end{remark}

    Theorem~\ref{Theorem:1.2} gives  sharp estimates to the life span $T(\kappa\psi)$ as $\kappa\to+0$ with $\psi$ behaving like $|x|^{-A} (A>0)$ at the space infinity. 
Appendix contains a brief summary of Theorem~\ref{Theorem:1.2} (see Tables 4 and 5 in Appendix).

\begin{theorem}
\label{Theorem:1.2}
Let $A>0$ and $\psi(x)=(1+|x|)^{-A}$.  
Then $T(\kappa\psi)\to\infty$ as $\kappa\to 0$ and the following holds: 
\begin{itemize}
  \item[{\rm (i)}] 
  Let $1<p<p_*$ or $0<A<1/(p-1)$. Then 
  $$
  T(\kappa\psi)\thicksim
  \left\{
  \begin{array}{ll}
  \kappa^{-\left(\frac{1}{2(p-1)}-\frac{1}{2} \min \{A,N\} \right)^{-1}} & \mbox{if}\quad A\not=N,\vspace{3pt}\\
  \displaystyle{\left(\frac{\kappa^{-1}}{\log(\kappa^{-1})}\right)^{\left(\frac{1}{2(p-1)}-\frac{N}{2}\right)^{-1}}} & 
  \mbox{if}\quad A=N,\\
  \end{array}
  \right.
  $$
  as $\kappa\to +0$;
  \item[{\rm (ii)}] 
  Let $p=p_*$ and $A\ge 1/(p-1)$. Then 
  $$
  \log T(\kappa\psi)\thicksim
  \left\{
  \begin{array}{ll}
  \kappa^{-(p-1)} & \mbox{if}\quad A>N,\vspace{3pt}\\
  \kappa^{-\frac{p-1}{p}} & \mbox{if}\quad A=N,\\
  \end{array}
  \right.
  $$
  as $\kappa\to +0$;
  \item[{\rm (iii)}] 
  Let $p>p_*$ and $A\ge 1/(p-1)$. 
  Then problem~\eqref{prob} with \eqref{ini}
  possesses a global-in-time solution if $\kappa>0$ is sufficiently small. 
\end{itemize}
\end{theorem}
\vspace{5pt}
\begin{remark}
Sharp estimates of the life span $T(\kappa\psi)$ as $\kappa\to+0$ have been already obtained in some cases. Specifically, if $\psi$ satisfies
$$
\psi(x)  = (1+|x|)^{-A} \quad (A>0)
$$
for all $x\in D$, then the following holds:
\begin{equation*}
T(\kappa\psi) \thicksim \left\{
\begin{array}{ll}
\kappa^{-\bigl(\frac{1}{2(p-1)} - \frac{A}{2}\bigr)^{-1}}\qquad & \mbox{if} \quad p\ge p_*,  \quad 0\le A <1/(p-1),\\
\kappa^{-\bigl(\frac{1}{2(p-1)} - \frac{1}{2} \min{\{A,\,N\}}\bigr)^{-1}}\qquad & \mbox{if} \quad p < p_*,   \quad A \neq N,\\
\displaystyle{\biggl(\frac{\kappa^{-1}}{\log(\kappa^{-1})}\biggr)^{\bigl(\frac{1}{2(p-1)} - \frac{N}{2}\bigr)^{-1}}}\qquad & \mbox{if} \quad p < p_*,   \quad A = N,\\
\end{array}
\right.
\qquad\qquad
\end{equation*}
as $\kappa\to+0$ \rm{(}see {\rm\cite{IS1}}\rm{)}.
\end{remark}

Finally, we show 
that $\lim_{\kappa\to 0}T(\kappa\psi)=\infty$ does not necessarily hold 
for problem~\eqref{prob} if $\psi$ has an exponential growth as $x_N\to\infty$. 
\begin{theorem}
\label{Theorem:1.3}
Let $p>1$, $\lambda>0$ and $\psi(x):=\exp{(\lambda x_N^2)}$.
Then 
\begin{equation}
\label{eq:1.6}
\lim_{\kappa\to+0} T(\kappa\psi)=(4\lambda)^{-1}.
\end{equation}
\end{theorem}

\begin{remark}
Let $\psi(x)=\exp{(\lambda x_N^2)}$. Set
$$
v(x,t):=\int_D G(x,y,t) \psi(y) \, dy.
$$ 
Then $v$ is a solution to 
\begin{equation*}
\left\{
\begin{array}{ll}
\partial_t v=\Delta v,\qquad & x\in{\bf R}^N_+,\,\,\,t>0,\vspace{3pt}\\
\displaystyle{-\frac{\partial v}{\partial x_N} }=0, & x\in\partial{\bf R}^N_+,\,\,\,t>0,\vspace{3pt}\\
v(x,0)=\psi(x),\qquad  & x\in \overline{{\bf R}^N_+},
\end{array}
\right.
\qquad\qquad
\end{equation*}
and 
$$
v(x,t) = (1-4\lambda t)^{-\frac{1}{2}} \exp\left(\frac{\lambda x_N^2}{1-4\lambda t}\right),
$$ 
where $N\ge1$. Moreover, $v$ does not exist after $t=(4\lambda)^{-1}$. 
\end{remark}

The rest of this paper is organized as follows. In Section 2, we review some of the facts on the solvability of  problem~\eqref{prob}, which have been already obtained in \cite{HI19}. In Section 3, we give  upper estimates and lower estimates to the life span $T(\kappa\psi)$ as $\kappa \to \infty$ (see Propositions~\ref{Proposition:3.1} and \ref{Proposition:3.2}). By combining these estimates, we can prove Theorem~\ref{Theorem:1.1}. In Section 4, we prove Theorem~\ref{Theorem:1.2} by the same method as in Section 3 (see Propositions~\ref{Proposition:4.1} and \ref{Proposition:4.2}) and  prove Theorem~\ref{Theorem:1.3}.
Appendix contains summaries of Theorems~\ref{Theorem:1.1} and \ref{Theorem:1.2}.

\section{Necessary conditions and sufficient conditions for the solvability of problem~\eqref{prob}}
In what follows the letter $C$ denotes a generic positive constant depending only on $N$ and $p$.
For any $L\ge0$, we set
\begin{eqnarray*}
& &
D_L:=\{(x',x_N) : x'\in{\bf R}^{N-1}, x_N\ge L^\half\},\\
& &
D_L':=\{(x',x_N) : x'\in{\bf R}^{N-1}, 0\le x_N <L^\half\}.
\end{eqnarray*}
Now we review  necessary conditions for the solvability of problem~\eqref{prob}, which have been obtained in \cite{HI19}.

\begin{theorem}
\label{Theorem:2.1}
Let $p>1$ and $u$ be a solution to \eqref{prob}  in $[0,T)$ with $u(0) = \varphi$, where $0<T<\infty$. Then for any $\delta>0$, there exists $\gamma_1 = \gamma_1(N,p,\delta)>0$ such that
\begin{equation}
\label{eq:2.1}
\sup_{x\in\RN} \exp\biggl(-(1+\delta)\frac{x_N^2}{4\sigma^2} \biggr) \int_{B_+(x,\sigma)} \varphi(y) \, dy \le \gamma_1 \sigma^{N-\frac{1}{p-1}}
\end{equation}
for $0<\sigma\le T^{1/2}$. In particular, in the case of $p=p_*$, there exists $\gamma_1' = \gamma_1'(N,\delta)>0$ such that
\begin{equation}
\label{eq:2.2}
\sup_{x\in\RN} \exp\biggl(-(1+\delta)\frac{x_N^2}{4\sigma^2} \biggr) \int_{B_+(x,\sigma)} \varphi(y) \, dy \le \gamma_1' \biggl[ \log\biggl( e + \frac{T^\half}{ \sigma}\biggr) \biggr]^{-N}
\end{equation}
for $0<\sigma\le T^{1/2}$.
\end{theorem}

\begin{remark}
If $1<p\le p_*$ and $\mu\not\equiv 0$ in $D$, 
then problem~\eqref{prob} possesses no nonnegative global-in-time solutions. 
See {\rm\cite{DFL}} and {\rm\cite{GL}}. 
\end{remark}
\vspace{5pt}

Next, we review  sufficient conditions for the solvability of problem~\eqref{prob}, which have been obtained also in \cite{HI19}. For any measurable function $\phi$ in $\RN$ and any bounded Borel set $E$, we set
$$
\dashint_E \phi(y) \, dy = \frac{1}{|E|} \int_E \phi(y) \, dy, \quad 
\phi_E(x) := \phi(x)\chi_{E}(x),
$$
where $|E|$ is the Lebesgue measure of $E$.

\begin{theorem}
\label{Theorem:2.2}
Let $1<p<p_*$, $T>0$ and $\delta\in(0,1)$. Set $\lambda := (1-\delta)/4T$. Then there exists $\gamma_2= \gamma_2(N,p,\delta)>0$ with the following property:
If $\varphi$ is a nonnegative measurable function in $\RNpls$ satisfying
\begin{equation}
\label{eq:2.3}
\sup_{x\in D} \dashint_{B_+(x,T^{1/2})} e^{-\lambda y_N^2} \varphi(y) \, dy \le \gamma_2 T^{-\frac{1}{2(p-1)}},
\end{equation}
then there exists a solution $u$ to \eqref{prob}  in $[0,T)$ with $u(0) = \varphi$. 
\end{theorem}
\vspace{5pt}

\begin{theorem}
\label{Theorem:2.3}
Let $p>1$, $a\in(1,p)$, $T>0$ and $\delta\in(0,1)$. Let $\varphi$ be a nonnegative measurable function in $\RNpls$. Set $\varphi_1 := \varphi_{D_T}$, $\varphi_2 := \varphi_{D_T'}$ and $\lambda := (1-\delta)/4T$. Then there exists $\gamma_3 = \gamma_3(N,p,a,\delta)>0$ with the following property:
Assume that $\varphi_1$ satisfies
\begin{equation}
\label{eq:2.4}
\sup_{x\in D} \dashint_{B_+(x,T^{1/2})} e^{-\lambda y_N^2} \varphi_1(y) \, dy \le \gamma_3 T^{-\frac{1}{2(p-1)}}.
\end{equation}
Furthermore, assume that $\varphi_2$ satisfies
\begin{equation}
\label{eq:2.5}
\sup_{x\in D_T'} \biggl[\dashint_{B_+(x,\sigma)} \varphi_2(y)^a \, dy \biggr]^\frac{1}{a} \le \gamma_3 \sigma^{-\frac{1}{p-1}} \quad \mbox{for } 0<\sigma\le T^\half.
\end{equation}
Then there exists a solution $u$ to \eqref{prob} in $[0,T)$ with $u(0) = \varphi$.
\end{theorem}
\vspace{5pt}

\begin{theorem}
\label{Theorem:2.4}
Let $p=p_*$, $T>0$ and $\delta\in(0,1)$. 
 Let $\varphi$ be a nonnegative measurable function in $\RNpls$. Set $\varphi_1 := \varphi_{D_T}$, $\varphi_2 := \varphi_{D_T'}$, $\lambda := (1-\delta)/4T$ and 
\begin{equation}
\label{eq:2.6}
\Phi(s):=s[\log (e+s)]^N,
\quad
\rho(s):=
s^{-N}\biggr[\log\biggr(e+\frac{1}{s}\biggr)\biggr]^{-N}
\quad\mbox{for}\quad s>0. 
\end{equation}
Then there exists $\gamma_4=\gamma_4(N,\delta)>0$ with the following property:
  Assume that $\varphi_1$ satisfies 
  \begin{equation}
  \label{eq:2.7}
  \sup_{x\in D_T}\,\dashint_{B_+(x,T^{1/2})}e^{-\lambda y_N^2}\varphi_1(y)\,dy
  \le\gamma_4T^{-\frac{1}{2(p-1)}}.
  \end{equation}
  Furthermore, assume that $\varphi_2$ satisfies 
  \begin{equation}
  \label{eq:2.8}
  \sup_{x\in D_T'}\Phi^{-1}\left[\,\dashint_{B_+(x,\sigma)}
  \Phi(T^\frac{1}{2(p-1)}\varphi_2(y))\,dy\,\right]\le\gamma_4\rho(\sigma T^{-\frac{1}{2}})
  \quad\mbox{for $0<\sigma\le T^{\frac{1}{2}}$}.
  \end{equation}
  Then there exists a solution $u$ to \eqref{prob} in $[0,T)$ with $u(0) = \varphi$.
\end{theorem}

\section{Proof of Theorem~\ref{Theorem:1.1}}
For simplicity of notation, we write $T_\kappa$ instead of  $T(\kappa\psi)$.
In this section we study the behavior of $T_\kappa$ as $\kappa\to\infty$ and prove Theorem~\ref{Theorem:1.1}. 
In order to prove Theorem~\ref{Theorem:1.1}, we obtain upper and lower estimates of $T_\kappa$ as $\kappa\to\infty$.
Proposition~\ref{Proposition:3.1} gives upper estimates of $T_\kappa$ as $\kappa\to\infty$.  
In the rest of this paper, 
for any two nonnegative functions
$f_1$ and $f_2$ defined in a subset $E$ of $[0,\infty)$,
we write
$f_1(t)\asymp f_2(t)$ for all $t\in E$ 
if $C^{-1}f_2(t)\le f_1(t)\le Cf_2(t)$ for all $t\in E$. 

\begin{proposition}
\label{Proposition:3.1}
Let $\psi$ be a nonnegative measurable function in $D$ such that 
\begin{equation}
\label{eq:3.1}
\psi(y)\ge |y|^{-A}\biggr[\log\left(e+\frac{1}{|y|}\right)\biggr]^{-B},
\qquad y\in B_+(0,1),
\end{equation}
where $0\le A\le N$ and $B$ are as in \eqref{eq:1.5}.
Then $\lim_{\kappa\to\infty}T(\kappa\psi)=0$. 
Furthermore, the following holds:
\vspace{3pt}
\newline
{\rm (i)} 
Let $1<p<p_*$. Then there exists $\gamma>0$ such that 
\begin{eqnarray}
\label{eq:3.2}
 & & 
T(\kappa\psi)\le\gamma [\kappa(\log\kappa)^{-B}]^{-\frac{2(p-1)}{-A(p-1)+1}}
\quad\qquad\mbox{if}\quad A<N,\,\, B\in{\bf R},\vspace{7pt}\\
\label{eq:3.3}
 & &
T(\kappa\psi)\le\gamma 
[\kappa(\log\kappa)^{-B+1}]^{-\frac{2(p-1)}{-A(p-1)+1}}
\qquad\mbox{if}\quad A=N,\,\, B>1,
\end{eqnarray}
for  sufficiently large $\kappa>0$;
\vspace{3pt}
\newline
{\rm (ii)} Let $p>p_*$. If, either 
\begin{equation}
\label{eq:3.4}
A>1/(p-1)  \quad \mbox{and} \quad B\in{\bf R} \qquad or \qquad A=1/(p-1) \quad \mbox{and} \quad B<0,
\end{equation}
then problem~\eqref{prob} with \eqref{ini} possesses no local-in-time solutions
for all $\kappa>0$. 
If 
$$
A=1/(p-1)\quad and \quad B=0, 
$$
then problem~\eqref{prob} with \eqref{ini} possesses no local-in-time solutions 
for sufficiently large $\kappa>0$. 
Furthermore, 
\begin{itemize}
  \item[{\rm (a)}] if $A<1/(p-1)$, then \eqref{eq:3.2} holds;
  \item[{\rm (b)}] if $A=1/(p-1)$ and $B>0$, then 
  there exists $\gamma'>0$ such that
  $$
  T(\kappa\psi)\le\exp(-\gamma'\kappa^{\frac{1}{B}})
  $$
  for  sufficiently large $\kappa>0$;
\end{itemize}
{\rm (iii)} Let $p=p_*$. If 
$$
A=N\quad \mbox{and}\quad B<N+1, 
$$
then problem~\eqref{prob} with \eqref{ini} possesses no local-in-time solutions  
for all $\kappa>0$. 
If 
$$
A=N\quad \mbox{and}\quad B=N+1, 
$$
then problem~\eqref{prob} with \eqref{ini} possesses no local-in-time solutions  
for  sufficiently large $\kappa>0$. Furthermore, 
\begin{itemize}
  \item[{\rm (c)}] if $A<N$, then \eqref{eq:3.2} holds;
  \item[{\rm (d)}] if $A=N$ and $B>N+1$, 
  then there exists $\gamma''>0$ such that
  $$
  T(\kappa\psi)\le\exp(-\gamma''\kappa^{\frac{1}{B-N-1}})
  $$
  for  sufficiently large $\kappa>0$. 
\end{itemize}
\end{proposition}
{\bf Proof.}
We assume that \eqref{prob} with \eqref{ini} possesses a solution in $[0,T_\kappa)$.
For any $p>1$, 
by \eqref{eq:2.1} and \eqref{eq:3.1}
we can find a constant $\gamma_1>0$ such that 
\begin{equation}
\label{eq:3.5}
\gamma_1\sigma^{N-\frac{1}{p-1}}\ge\kappa\int_{B_+(0,\sigma)}\psi(y)\,dy
\ge \kappa\int_{B_+(0,\sigma)}|y|^{-A}\biggr[\log\left(e+\frac{1}{|y|}\right)\biggr]^{-B}\,dy>0
\end{equation}
for $0<\sigma\le T_\kappa^{1/2}$. 
Firstly, we show that 
$\lim_{\kappa\to\infty}T_\kappa=0$ 
by contradiction.
Assume that there exist $\{\kappa_j\}_{j=1}^\infty$ and $c_*>0$ such that
$$
\lim_{j\to\infty} \kappa_j = \infty, \quad T_{\kappa_j} > c_*^2 \quad \mbox{for all } j=1,2,\cdots.
$$
By \eqref{eq:3.5} with $\sigma=c_*$, we have
$$
\gamma_1c_*^{N-\frac{1}{p-1}} \ge \kappa_j\int_{B_+(0,c_*)}|y|^{-A}\biggr[\log\left(e+\frac{1}{|y|}\right)\biggr]^{-B}\,dy, \quad j=1,2,\cdots,
$$
where $\gamma_1$ is a constant independent of $\kappa_j$. 
Since $\lim_{j\to\infty} \kappa_j =\infty$, we have a contradiction.
Since $c_*$ is arbitrary, we have
$$
\lim_{\kappa\to\infty} T_\kappa = 0.
$$
Without loss of generality we can assume that $T_\kappa>0$ is sufficiently small.

We prove assertion~(i). 
Let $1<p<p_*$. 
For any $p>1$, by \eqref{eq:3.5} we have 
\begin{equation}
\label{eq:3.6}
\gamma_1\ge 
\left\{
\begin{array}{ll}
C\kappa \sigma^{-A+\frac{1}{p-1}}\biggr[\log\left(e+\sigma^{-1}\right)\biggr]^{-B} & \mbox{if}\quad A<N,\,\,B\in{\R},\vspace{7pt}\\
C\kappa\sigma^{-N+\frac{1}{p-1}} \biggr[\log\left(e+\sigma^{-1}\right)\biggr]^{-B+1} & \mbox{if}\quad A=N,\,\,B>1,
\end{array}
\right.
\end{equation}
for $0<\sigma\le T_\kappa^{1/2}$ and sufficiently large $\kappa>0$. 
We notice that for any $a_1>0$ and $a_2\in\R$,
\begin{equation}
\label{eq:3.7}
\Psi(\tau) := \tau^{a_1} [\log(e+\tau^{-1})]^{a_2} 
\,\,\mbox{is increasing for sufficiently small}\,\, \tau>0,
\end{equation}
$\Psi^{-1}$ satisfies
\begin{equation}
\label{eq:3.8}
\Psi^{-1}(\tau)\asymp\tau^\frac{1}{a_1} [\log(e+\tau^{-1})]^{-\frac{a_2}{a_1}} \quad \rm{for\,\,sufficiently \,\, small} \,\, \tau>0
\end{equation}
 and $\Psi^{-1}(\tau)$ is also increasing sufficiently small $\tau>0$.
We consider the case where $A<N$ and $B\in{\bf R}$. Set
$$
a_1:=-A+\frac{1}{p-1}>0 \quad \mbox{and} \quad a_2:= -B.
$$
By \eqref{eq:3.6}, \eqref{eq:3.7} and \eqref{eq:3.8} we have
\begin{equation*}
\begin{split}
\sigma 
&\le C\Psi^{-1}(C\gamma_1\kappa^{-1})\\
&\le C(C\gamma_1\kappa^{-1})^{(-A+\frac{1}{p-1})^{-1}}[\log(e+(C\gamma_1\kappa^{-1})^{-1})]^{B(-A+\frac{1}{p-1})^{-1}}\\
&\le C[\kappa(\log\kappa)^{-B}]^{-\frac{p-1}{-A(p-1)+1}}
\end{split}
\end{equation*}
for $0<\sigma\le T_\kappa^{1/2}$. Setting $\sigma=T_\kappa^{1/2}$, we obtain \eqref{eq:3.2}. Similarly, we can obtain \eqref{eq:3.3}.
Thus assertion~(i) follows. 

We prove assertion~(ii). Let $p>p_*$. 
We can assume that
\begin{equation}
\label{eq:3.9}
A<N\quad\mbox{and}\quad B\in\R \qquad \mbox{or} \qquad 
A=N \quad\mbox{and}\quad B>1
\end{equation}
because
$$
\int_{B_+(0,\sigma)} |y|^{-A}\biggr[\log\left(e+\frac{1}{|y|}\right)\biggr]^{-B}\,dy=\infty
$$
for all $\sigma>0$ when $A$ and $B$ do not satisfy \eqref{eq:3.9}. 
By \eqref{eq:3.5}, this implies that $T_\kappa=0$ for all $\kappa>0$.
By \eqref{eq:3.9}, we have \eqref{eq:3.6}.
Since $A$ and $B$ satisfy \eqref{eq:3.4}, the right hand side of \eqref{eq:3.6} goes to infinity as $\sigma\to+0$. This implies that $T_\kappa=0$ for all $\kappa>0$.
In the case where $A=1/(p-1)$ and $B=0$ (this condition also satisfies \eqref{eq:3.9}), it follows from \eqref{eq:3.6} that 
\begin{equation}
\label{eq:3.9.1}
\gamma_1\ge C\kappa.
\end{equation}
Since \eqref{eq:3.9.1} does not hold for sufficiently large $\kappa>0$, this implies that $T_\kappa=0$ for  sufficiently large $\kappa>0$. 
Furthermore, if $A<1/(p-1)$, we obtain \eqref{eq:3.2} by a similar argument to the proof of assertion~(i). 
Then we obtain (a).
It remains to consider the case where $A=1/(p-1)$ and $B>0$.
Since $T_\kappa>0$ is sufficiently small, by \eqref{eq:3.6} we have
$$
\gamma\kappa^{-1}\ge C[\log(e+T_\kappa^{-\frac{1}{2}})]^{-B}\ge C[\log(T_\kappa^{-\frac{1}{2}})]^{-B}.
$$
 Since $B>0$, this implies that there  exists a constant $\gamma'>0$ such that
$$
T_\kappa\le \exp(-\gamma'\kappa^\frac{1}{B})
$$
for sufficiently large $\kappa>0$ and (b) follows. Thus assertion~(ii) is proved.

Finally, we prove assertion~(iii).
Let $A=N$. 
Since $p=p_*$ and $B>1$, by \eqref{eq:2.2} we have
\begin{equation}
\label{eq:3.10}
\begin{split}
\gamma_1'\biggr[\log\biggr(e+\frac{T_\kappa^{\frac{1}{2}}}{\sigma}\biggr)\biggr]^{-N}
 & \ge\kappa\int_{B_+(0,\sigma)}\psi(y)\,dy\\
& \ge \kappa\int_{B_+(0,\sigma)}|y|^{-A}\biggr[\log\left(e+\frac{1}{|y|}\right)\biggr]^{-B}\,dy\\
 & \ge C\kappa[\log(e+\sigma^{-1})]^{-B+1}
\end{split}
\end{equation}
for $0<\sigma\le T_\kappa^{1/2}$.
In the case of $B<N+1$, we see that \eqref{eq:3.10} does not hold for sufficiently small $\sigma>0$. 
This implies that $T_\kappa=0$ for all $\kappa>0$. 
In the case of $B=N+1$, it follows from \eqref{eq:3.10} with $\sigma=T_\kappa (<T_\kappa^{1/2})$ that
$$
\gamma_1'[\log(e+T_\kappa^{-\frac{1}{2}})]^{-N}\ge C\kappa[\log(e+T_\kappa^{-1})]^{-N}.
$$
This inequality implies that 
\begin{equation}
\label{eq:3.10.1}
\gamma_1'\ge C \kappa.
\end{equation}
Since \eqref{eq:3.10.1} does not hold for  sufficiently large $\kappa>0$, this implies that $T_\kappa=0$ for  sufficiently large $\kappa>0$. 
In the case of $A<N$, since \eqref{eq:3.6} holds, we obtain \eqref{eq:3.2} by a similar argument to the proof of assertion~(i).  
Then we obtain (c).
In the case where $A=N$ and $B>N+1$, 
since $T_\kappa>0$ is sufficiently small, by \eqref{eq:3.10} with $\sigma=T_\kappa(<T_\kappa^{1/2})$ we have
$$
C\gamma_1'\kappa^{-1}\le [\log(e+T_\kappa^{-1})]^{-B+N+1} \le [\log(T_\kappa^{-1})]^{-B+N+1}.
$$
Since $B-N-1>0$, this implies that there exists a constant $\gamma''>0$ such that 
$$
T_\kappa\le \exp(-\gamma''\kappa^\frac{1}{B-N-1})
$$
for sufficiently large $\kappa>0$.
Thus assertion~(iii) follows and the proof of Proposition~\ref{Proposition:3.1} is complete. 
$\Box$
\vspace{5pt}

In Proposition~\ref{Proposition:3.2} 
we obtain lower estimates of $T_\kappa$ as $\kappa\to\infty$ 
and show the optimality of the estimates of $T_\kappa$ 
in Proposition~\ref{Proposition:3.1}. 

\begin{proposition}
\label{Proposition:3.2}
Let $\psi$ be a nontirivial nonnegative measurable function in $D$ such that 
$\mbox{supp}\,\psi\subset B(0,1)$ and 
\begin{equation}
\label{eq:3.11}
\psi(y)\le |y|^{-A}\biggr[\log\left(e+\frac{1}{|y|}\right)\biggr]^{-B},
\qquad y\in B_+(0,1),
\end{equation}
where $0\le A\le N$ and $B$ are as in \eqref{eq:1.5}. 
\vspace{3pt}
\newline
{\rm (i)} 
Let $1<p<p_*$. Then there exists $\gamma>0$ such that 
\begin{eqnarray}
\label{eq:3.12}
 & & 
T(\kappa\psi) \ge\gamma[\kappa(\log\kappa)^{-B}]^{-\frac{2(p-1)}{-A(p-1)+1}}
\quad\qquad\mbox{if}\quad A<N,\,\, B\in{\bf R},\vspace{7pt}\\
\nonumber 
 & & 
T(\kappa\psi) \ge \gamma[\kappa(\log\kappa)^{-B+1}]^{-\frac{2(p-1)}{-A(p-1)+1}}
\qquad\mbox{if}\quad A=N,\,\, B>1.
\end{eqnarray}
{\rm (ii)} Let $p>p_*$. 
\begin{itemize}
  \item[{\rm (a)}] If $A<1/(p-1)$, then \eqref{eq:3.12} holds;
  \item[{\rm (b)}] If $A=1/(p-1)$ and $B>0$, then there exists $\gamma'>0$ such that
  $$
  T(\kappa\psi)\ge\exp(-\gamma'\kappa^{\frac{1}{B}})
  $$
  for  sufficiently large $\kappa>0$; 
\end{itemize}
{\rm (iii)} Let $p=p_*$. 
\begin{itemize}
  \item[{\rm (c)}] If $A<N$, then \eqref{eq:3.12} holds;
  \item[{\rm (d)}] If $A=N$ and $B>N+1$, 
  then there exists $\gamma''>0$ such that
  $$
  T(\kappa\psi)\ge\exp(-\gamma''\kappa^{\frac{1}{B-N-1}})
  $$
  for sufficiently large $\kappa>0$. 
  \end{itemize}
\end{proposition}
{\bf Proof.}
We first consider the case where $p>p_*$ and $A<1/(p-1)(<N)$. 
Let $a\in(1,p)$ be such that $a A<N$. 
By the Jensen inequality and \eqref{eq:3.11}, we have
\begin{equation}
\label{eq:3.13}
\begin{split}
  \sigma^{\frac{1}{p-1}} \sup_{x\in D}\,\dashint_{B_+(x,\sigma)}
\kappa\psi(y)\,dy
&\le \sigma^{\frac{1}{p-1}}\sup_{x\in D}\left[\,\dashint_{B_+(x,\sigma)}
[\kappa\psi(y)]^a\,dy\,\right]^{\frac{1}{a}}\\
 & \le C\kappa\sigma^{\frac{1}{p-1}}\biggr[\,\dashint_{B_+(0,\sigma)}|y|^{-Aa}
\biggr[\log\left(L+\frac{1}{|y|}\right)\biggr]^{-a B}\,dy\biggr]^{\frac{1}{a}}\\
&\le C\kappa\sigma^{\frac{1}{p-1}-A}\biggr[\log\left(e+\frac{1}{\sigma}\right)\biggr]^{-B}
\end{split}
\end{equation}
for  sufficiently small $\sigma>0$. 
Let $c$ be a sufficiently small positive constant and set  
$$
\tilde{T}_\kappa:=c[\kappa(\log\kappa)^{-B}]^{-\frac{2(p-1)}{-A(p-1)+1}}.
$$
Since $A<1/(p-1)$, by \eqref{eq:3.7} we have
\begin{equation*}
\begin{split}
C\kappa\sigma^{\frac{1}{p-1}-A}\biggr[\log\left(e+\frac{1}{\sigma}\right)\biggr]^{-B}
&\le C\kappa\sigma^{\frac{1}{p-1}-A}\biggr[\log\left(e+\frac{1}{\sigma}\right)\biggr]^{-B}\biggr|_{\sigma=\tilde{T}_\kappa^{1/2}}\\
&\le C c^{\frac{1}{2(p-1)} -\frac{A}{2}} 
\end{split}
\end{equation*}
for $0<\sigma\le\tilde{T}_\kappa^{1/2}$ and  sufficiently large $\kappa>0$.
Taking a sufficiently small $c>0$ if necessary, we obtain
\begin{equation}
\label{eq:3.14}
C\kappa\sigma^{\frac{1}{p-1}-A}\biggr[\log\left(e+\frac{1}{\sigma}\right)\biggr]^{-B}\le \gamma_3,
\end{equation} 
for $0<\sigma\le\tilde{T}_\kappa^{1/2}$ and  sufficiently large $\kappa>0$, where $\gamma_3$ is as in Theorem~\ref{Theorem:2.3}. 
Then \eqref{eq:3.13} and \eqref{eq:3.14} yield \eqref{eq:2.4} and \eqref{eq:2.5}.  Applying Theorem~\ref{Theorem:2.3}, we see that 
\eqref{prob} with \eqref{ini} possesses a solution in $[0,\tilde{T}_\kappa )$ and
$$
T_\kappa\ge\tilde{T}_\kappa = c[\kappa(\log\kappa)^{-B}]^{-\frac{2(p-1)}{-A(p-1)+1}}
$$ 
for  sufficiently large $\kappa>0$. 
So we have (a). Similarly, we have (b) and (c). 
Furthermore, 
we can prove \eqref{eq:3.12} in the case of $1<p<p_*$ by using the above argument with $a=1$ 
and applying Theorem~\ref{Theorem:2.2}. 

Next we consider the case where $1<p<p_*$, $A=N$ and $B>1$, let $c$ be a sufficiently small positive constant and set
$$
\tilde{T}_\kappa':=c[\kappa(\log\kappa)^{-B+1}]^{-\frac{2(p-1)}{-A(p-1)+1}}.
$$
Taking a sufficiently small $c>0$ if necessary, by \eqref{eq:3.11} we have
\begin{equation}
\label{eq:3.15}
\begin{split}
 & \tilde{T}_\kappa'^{\frac{1}{2(p-1)}} \sup_{x\in D}\,\dashint_{B_+(x,\tilde{T}_\kappa'^{1/2})}
\kappa\psi(y)\,dy\\
&\le C\kappa\tilde{T}_\kappa'^{\frac{1}{2(p-1)}}\dashint_{B_+(0,\tilde{T}_\kappa'^{1/2})}
|y|^{-N}\biggr[\log\left(\frac{1}{|y|}\right)\biggr]^{-B}\,dy\\
 & \le C\kappa\tilde{T}_\kappa'^{\frac{1}{2(p-1)}-\frac{N}{2}}\biggr[\log\left(\frac{1}{\tilde{T}_\kappa'^{1/2}}\right)\biggr]^{-B+1}
 \le Cc^{\frac{1}{2(p-1)}-\frac{N}{2}}\le\gamma_2
\end{split}
\end{equation}
for  sufficiently large $\kappa>0$, where $\gamma_2$ is as in Theorem~\ref{Theorem:2.2}. Then \eqref{eq:3.15} yields \eqref{eq:2.3}.
Applying Theorem~\ref{Theorem:2.2}, we see that \eqref{prob} with \eqref{ini} possesses a solution in $[0,\tilde{T}_\kappa')$ and
$$
T_\kappa\ge \tilde{T}_\kappa'=c[\kappa(\log\kappa)^{-B+1}]^{-\frac{2(p-1)}{-A(p-1)+1}}
$$
 for sufficiently large $\kappa>0$. So we have assertion~(i).
 
 It remains to prove (d). 
Let $p=p_*$, $A=N$ and $B>N+1$. 
Let $c$ be a sufficiently small positive constant and 
set 
$$
\hat{T}_\kappa:=\exp(-c^{-1}\kappa^{\frac{1}{B-N-1}})
$$ 
for sufficiently large $\kappa>0$.
We can assume that $\hat{T}_\kappa>0$ is sufficiently small.
Let $\Phi$ be as in \eqref{eq:2.6}.
 By \eqref{eq:3.7} and \eqref{eq:3.8} with $a_1=1$ and $a_2=N$, we see that $\Phi^{-1}$ satisfies
$$
\Phi^{-1}(\tau)\asymp\tau [\log(e+\tau^{-1})]^{-N} \quad \rm{for \,\, sufficiently \,\, small} \,\,\, \tau>0
$$
 and $\Phi^{-1}(\tau)$ is  increasing for sufficiently small $\tau>0$.
Similarly to \eqref{eq:3.13},  we have
\begin{equation}
\label{eq:3.16}
\begin{split}
 & \sup_{x\in D}\Phi^{-1}\left[\,\dashint_{B_+(x,\sigma)}
\Phi\left(\hat{T}_\kappa^{\frac{1}{2(p-1)}}\kappa\psi(y)\right)\,dy\right]\\
 & \qquad\quad
\le\Phi^{-1}\left[\,\dashint_{B_+(0,\sigma)}
\Phi\left(\hat{T}_\kappa^{\frac{N}{2}}\kappa|y|^{-N}\biggr[\log\biggr(e+\frac{1}{|y|}\biggr)\biggr]^{-B}\right)\,dy\right]
\end{split}
\end{equation}
for $0<\sigma\le \hat{T}_\kappa^{1/2}$.
Since 
\begin{equation*}
\begin{split}
 & \log\biggr[e+\hat{T}_\kappa^{\frac{N}{2}}\kappa|y|^{-N}\biggr[\log\biggr(e+\frac{1}{|y|}\biggr)\biggr]^{-B}\biggr]\\
 & \le\log\biggr[\biggr(e+\hat{T}_\kappa^{\frac{N}{2}}\kappa|y|^{-N}\biggr)\biggr(e+\biggr[\log\biggr(e+\frac{1}{|y|}\biggr)\biggr]^{-B}\biggr)\biggr]\\
 & \le\log\biggr[C\hat{T}_\kappa^{\frac{N}{2}}\kappa|y|^{-N}\biggr]\le C\log\biggr[\hat{T}_\kappa^{\frac{N}{2}}\kappa|y|^{-N}\biggr]
 \le C\log\frac{1}{|y|}
\end{split}
\end{equation*}
for $y\in B_+(0,\sigma)$, $0<\sigma\le\hat{T}_\kappa^{1/2}$ and  sufficiently large $\kappa$, 
we have 
\begin{equation*}
\begin{split}
 & \dashint_{B_+(0,\sigma)}\Phi\left(\hat{T}_\kappa^{\frac{N}{2}}\kappa|y|^{-N}\biggr[\log\biggr(e+\frac{1}{|y|}\biggr)\biggr]^{-B}\right)\,dy\\
 &= 	\dashint_{B_+(0,\sigma)}\hat{T}_\kappa^{\frac{N}{2}}\kappa|y|^{-N}\biggr[\log\biggr(e+\frac{1}{|y|}\biggr)\biggr]^{-B} \left[\log\biggr[e+\hat{T}_\kappa^{\frac{N}{2}}\kappa|y|^{-N}\biggr[\log\biggr(e+\frac{1}{|y|}\biggr)\biggr]^{-B}\biggr]\right]^N\,dy\\
 & \le 	C\hat{T}_\kappa^{\frac{N}{2}}\kappa\,\,
\dashint_{B_+(x,\sigma)}\,|y|^{-N}\biggr[\log\frac{1}{|y|}\biggr]^{-B+N}\,dy
\le C\kappa\sigma^{-N}\hat{T}_\kappa^{\frac{1}{2(p-1)}}\biggr[\log\frac{1}{\sigma}\biggr]^{-B+N+1}
\end{split}
\end{equation*}
for $0<\sigma\le\hat{T}_\kappa^{1/2}$ and  sufficiently large $\kappa>0$. 
This together with \eqref{eq:3.16} implies that 
\begin{equation}
\label{eq:3.17}
\begin{split}
 & \sup_{x\in D}\dashint_{B_+(x,\sigma)}\hat{T}_\kappa^{\frac{1}{2(p-1)}}\kappa\psi(y)\,dy
\le\sup_{x\in D}\Phi^{-1}\left[\,\dashint_{B_+(x,\sigma)}
\Phi\left(\hat{T}_\kappa^{\frac{1}{2(p-1)}}\kappa\psi(y)\right)\,dy\right]\\
 & \le \Phi^{-1}\left(\dashint_{B_+(0,\sigma)}\Phi\left(\hat{T}_\kappa^{\frac{N}{2}}\kappa|y|^{-N}\biggr[\log\biggr(e+\frac{1}{|y|}\biggr)\biggr]^{-B}\right)\,dy\right)\\
 & \le \Phi^{-1}\left( C\kappa\sigma^{-N}\hat{T}_\kappa^{\frac{1}{2(p-1)}}\biggr[\log\frac{1}{\sigma}\biggr]^{-B+N+1}\right)\\
 & \le C\kappa\sigma^{-N}\hat{T}_\kappa^{\frac{N}{2}}\biggr[\log\frac{1}{\sigma}\biggr]^{-B+N+1}
 \biggr(\log\biggr[e+C\hat{T}_\kappa^{\frac{N}{2}}\kappa\sigma^{-N}\biggr[\log\frac{1}{\sigma}\biggr]^{-B+N+1}\biggr]\biggr)^{-N}\\
 & \le C\kappa\sigma^{-N}\hat{T}_\kappa^{\frac{N}{2}}\biggr[\log\frac{1}{\sigma}\biggr]^{-B+1}
\end{split}
\end{equation}
for $0<\sigma\le\hat{T}_\kappa^{1/2}$ and sufficiently large $\kappa>0$. 
On the other hand, since $\hat{T}_\kappa>0$ is sufficiently small, we have
\begin{equation}
\label{eq:3.18}
\rho(\sigma\hat{T}_\kappa^{-\frac{1}{2}})=\sigma^{-N}\hat{T}_\kappa^{\frac{N}{2}}
\biggr[\log\biggr(e+\frac{\hat{T}_\kappa^{\frac{1}{2}}}{\sigma}\biggr)\biggr]^{-N}
\ge\sigma^{-N}\hat{T}_\kappa^{\frac{N}{2}}\biggr[\log\frac{1}{\sigma}\biggr]^{-N}
\end{equation}
for $0<\sigma\le\hat{T}_\kappa^{1/2}$ and sufficiently large $\kappa$, where $\rho$ is as in \eqref{eq:2.6}. 
Since $B>N+1$ and
\begin{equation}
\label{eq:3.19}
\kappa\biggr[\log\frac{1}{\sigma}\biggr]^{-B+1+N}\le C\kappa\biggr[\log\frac{1}{\hat{T}_\kappa^{\frac{1}{2}}}\biggr]^{-B+1+N}
=Cc^{B-N-1},
\end{equation}
taking a sufficiently small $c>0$ if necessary, \eqref{eq:3.17}, \eqref{eq:3.18} and \eqref{eq:3.19} yield \eqref{eq:2.7} and \eqref{eq:2.8}.
Applying Theorem~\ref{Theorem:2.4}, we see that 
$$
T_\kappa\ge\hat{T}_\kappa = \exp(-c^{-1}\kappa^{\frac{1}{B-N-1}})
$$ 
for sufficiently large $\kappa>0$. 
This implies (d). The proof of Proposition~\ref{Proposition:3.2} is complete. 
$\Box$
\vspace{5pt}

\section{Proofs of Theorem~\ref{Theorem:1.2} and Theorem~\ref{Theorem:1.3}}
We state two results on the behavior of $T_\kappa$ as $\kappa\to +0$. 
If $\psi$ is a bounded function in ${\bf R}^N$, then $T_\kappa\to\infty$ as $\kappa\to +0$  
and the behavior of $T_\kappa$ depends on the decay of $\psi$ at the space infinity.  
In order to prove Theorem~\ref{Theorem:1.2}, It suffice to prove the following propositions.
In Proposition~\ref{Proposition:4.1} we obtain  upper estimates of $T_\kappa$ as $\kappa\to +0$.
\begin{proposition}
\label{Proposition:4.1}
Let $N\ge 1$ and $p>1$. 
Let $A>0$ and $\psi$ be a nonnegative $L^\infty(D)$-function such that 
$\psi(x)\ge(1+|x|)^{-A}$ for $x\in D$.
\vspace{3pt}
\newline
{\rm (i)}
Let $p=p_*$ and $A\ge 1/(p-1)=N$. 
Then there exists $\gamma>0$ such that 
$$
\log T(\kappa\psi)\le
\left\{
\begin{array}{ll}
  \gamma\kappa^{-(p-1)} & \mbox{if}\quad A>N,\vspace{3pt}\\
  \gamma\kappa^{-\frac{p-1}{p}} & \mbox{if}\quad A=N,\\
\end{array}
\right.
$$
for  sufficiently small $\kappa>0$. 
\vspace{3pt}
\newline
{\rm (ii)} 
Let $1<p<p_*$ or  $A<1/(p-1)$. 
Then there exists $\gamma'>0$ such that 
$$
T(\kappa\psi)\le
\left\{
\begin{array}{ll}
  \gamma'\kappa^{-\left(\frac{1}{2(p-1)}-\frac{1}{2} \min \{A,N\} \right)^{-1}} & \mbox{if}\quad A\not=N,\vspace{3pt}\\
  \displaystyle{\gamma'\left(\frac{\kappa^{-1}}{\log(\kappa^{-1})}\right)^{\left(\frac{1}{2(p-1)}-\frac{N}{2}\right)^{-1}}} & \mbox{if}\quad A=N,\\
\end{array}
\right.
$$
for  sufficiently small $\kappa>0$. 
\end{proposition}
{\bf Proof.}
Since $\psi\in L^\infty(D)$, by Theorem~\ref{Theorem:2.3} we have
$$
T_\kappa\ge C\kappa^{-(p-1)}
$$
for sufficiently small $\kappa>0$. This implies that $\lim_{\kappa\to0} T_\kappa = \infty$.
Without loss of generality, we can assume that $T_\kappa>0$ is sufficiently large.
For any $p>1$, we see that 
\begin{equation}
\begin{split}
\label{eq:4.1}
\int_{B_+(0,\sigma)}\kappa\psi(y)\,dy 
	&\ge\kappa\int_{B_+(0,\sigma)}(1+|y|)^{-A}\,dy\\
	&\ge\left\{
\begin{array}{ll}
  C\kappa				& \mbox{if}\quad \sigma>1, A>N,\vspace{3pt}\\
  C\kappa\log(e+\sigma)	& \mbox{if}\quad \sigma>1, A=N,\vspace{3pt}\\
  C\kappa\sigma^{N-A}	& \mbox{if}\quad \sigma>1, A<N,
  \end{array}
\right.
\end{split}
\end{equation}
for $\sigma>1$ and sufficiently small $\kappa>0$. 
In the case of $p=p_*$, it follows from \eqref{eq:2.2} that 
\begin{equation*}
\int_{B_+(0,\sigma)} \kappa\psi(y)\,dy \le \gamma_1'\left[\log\biggl(e+\frac{T_\kappa^\half}{\sigma}\biggr)\right]^{-N}
\end{equation*}
for $0<\sigma\le T_\kappa^{1/2}$ and sufficiently small $\kappa>0$. This implies that
\begin{equation}
\label{eq:4.2}
\int_{B_+(0,T_\kappa^{1/4})} \kappa\psi(y)\,dy \le C\gamma_1'[\log T_\kappa]^{-N},
\end{equation}
\begin{equation}
\label{eq:4.3}
\int_{B_+(0,T_\kappa^{1/2})} \kappa\psi(y)\,dy \le C\gamma_1',
\end{equation}
for sufficiently small $\kappa>0$.
By \eqref{eq:4.1} and \eqref{eq:4.2} with $\sigma=T_\kappa^{1/4}$ we obtain assertion~(i).
Furthermore, by \eqref{eq:4.1} and \eqref{eq:4.3} with $\sigma=T_\kappa^{1/2}$ we obtain assertion~(ii) in the case where $p=p_*$ and $A<1/(p-1)$. 

We prove assertion~(ii) in the case of $1<p<p_*$.
By \eqref{eq:2.1}
we see that 
\begin{equation}
\label{eq:4.4}
\int_{B_+(0,T_\kappa^{1/2})}\kappa\psi(y)\,dy
\le \gamma_1T_\kappa^{\frac{N}{2}-\frac{1}{2(p-1)}}. 
\end{equation}
By \eqref{eq:4.1} and \eqref{eq:4.4}, we obtain assertion~(ii) in the case of $1<p<p_*$. 
Similarly, we obtain assertion~(ii) in the case of $p>p_*$. Thus Proposition~\ref{Proposition:4.1} follows.
$\Box$
\vspace{5pt}

In Proposition~\ref{Proposition:4.2} 
we obtain lower estimates of $T_\kappa$ as $\kappa\to+ 0$ 
and show the optimality of the estimates of $T_\kappa$ 
in Proposition~\ref{Proposition:4.1}. 
\begin{proposition}
\label{Proposition:4.2}
Let $N\ge1$ and $p>1$. 
Let $A>0$ and $\psi$ be a nonnegative measurable function in $D$ such that 
$\mbox{supp}\,\psi\subset D$ and 
$0\le\psi(x)\le(1+|x|)^{-A}$ for $x\in D$.
\vspace{3pt}
\newline
{\rm (i)} 
Let $p=p_*$ and $A\ge 1/(p-1)=N$. 
Then there exists $\gamma>0$ such that 
$$
\log T(\kappa\psi)\ge
\left\{
\begin{array}{ll}
  \gamma\kappa^{-(p-1)} & \mbox{if}\quad A>N,\vspace{3pt}\\
  \gamma\kappa^{-\frac{p-1}{p}} & \mbox{if}\quad A=N,\\
\end{array}
\right.
$$
for sufficiently small $\kappa>0$. 
\vspace{3pt}
\newline
{\rm (ii)}
Let $1<p<p_*$ or  $A<1/(p-1)$. 
Then there exists $\gamma'>0$ such that 
$$
T(\kappa\psi)\ge
\left\{
\begin{array}{ll}
  \gamma'\kappa^{-\left(\frac{1}{2(p-1)}-\frac{1}{2} \min \{A,N\} \right)^{-1}} & \mbox{if}\quad A\not=N,\vspace{3pt}\\
  \displaystyle{\gamma'\left(\frac{\kappa^{-1}}{\log(\kappa^{-1})}\right)^{\left(\frac{1}{2(p-1)}-\frac{N}{2}\right)^{-1}}} & \mbox{if}\quad A=N,\\ \end{array}
\right.
$$
for sufficiently small $\kappa>0$. 
\end{proposition}
{\bf Proof.} 
Let $p=p_*$ and $A>N$. Let $c$ be a sufficiently small positive constant and set 
$$
\hat{T}_\kappa:=\exp(c\kappa^{-(p-1)})=\exp(c\kappa^{-\frac{1}{N}}).
$$
Let $L\ge e$ be such that
$$
\tau[\log(L+\tau)]^{-N} \mbox{is increasing in } [0,\infty).
$$
Then we see that $\Phi(\tau)\asymp \tau[\log(L+\tau)]^{N}$ and $\Phi^{-1}(\tau) \asymp \tau[\log(e+\tau)]^{-N} \asymp \tau[\log(L+\tau)]^{-N}$ for all $\tau>0$. 
Similarly to \eqref{eq:3.17}, we have
\begin{equation}
\label{eq:4.5}
\begin{split}
\sup_{x\in D}\dashint_{B_+(x,\sigma)}\hat{T}_\kappa^{\frac{1}{2(p-1)}}\kappa\psi(y)\,dy
 & \le\sup_{x\in D}\Phi^{-1}\left[\,\dashint_{B_+(x,\sigma)}
\Phi\left(\hat{T}_\kappa^{\frac{1}{2(p-1)}}\kappa\psi(y)\right)\,dy\right]\\
 & 
\le\Phi^{-1}\left[\,\dashint_{B_+(0,\sigma)}
\Phi\left(\hat{T}_\kappa^{\frac{N}{2}}\kappa (1+|y|)^{-A}\right)\,dy\right]
\end{split}
\end{equation}
for all $\sigma>0$. 
Since
\begin{equation}
\label{eq:4.6}
\log\biggr[L+\hat{T}_\kappa^{\frac{N}{2}}\kappa(1+|y|)^{-A}\biggr]
\le \log(C\hat{T}_\kappa^{\frac{N}{2}})\le Cc\kappa^{-\frac{1}{N}}
\end{equation}
for sufficiently small $\kappa>0$, 
we have 
\begin{equation}
\label{eq:4.7}
\dashint_{B_+(0,\sigma)}\Phi\left(\hat{T}_\kappa^{\frac{N}{2}}\kappa(1+|y|)^{-A}\right)\,dy\le Cc^N\hat{T}_\kappa^{\frac{N}{2}}\,\,
\dashint_{B_+(0,\sigma)}\,(1+|y|)^{-A}\,dy
\le Cc^N\hat{T}_\kappa^{\frac{N}{2}}\sigma^{-N}
\end{equation}
for $0<\sigma\le\hat{T}_\kappa^{1/2}$ and  sufficiently small $\kappa>0$. 
This together with \eqref{eq:4.5} implies that 
\begin{equation*}
\begin{split}
 & \sup_{x\in D}\dashint_{B_+(x,\sigma)}\hat{T}_\kappa^{\frac{1}{2(p-1)}}\kappa\psi(y)\,dy
\le\sup_{x\in D}\Phi^{-1}\left[\,\dashint_{B_+(x,\sigma)}
\Phi\left(\hat{T}_\kappa^{\frac{1}{2(p-1)}}\kappa\psi(y)\right)\,dy\right]\\
 & \le Cc^N\sigma^{-N}\hat{T}_\kappa^{\frac{N}{2}}
 \biggr(\log\biggr[L+Cc^N\hat{T}_\kappa^{\frac{N}{2}}\sigma^{-N}\biggr]\biggr)^{-N}
 \le Cc^N\sigma^{-N}\hat{T}_\kappa^{\frac{N}{2}}
 \biggr(\log\biggr[L+\frac{\hat{T}_\kappa^{\frac{1}{2}}}{\sigma}\biggr]\biggr)^{-N}\\
 & =Cc^N\rho(\sigma\hat{T}_\kappa^{-\frac{1}{2}})
\end{split}
\end{equation*}
for $0<\sigma\le\hat{T}_\kappa^{1/2}$ and sufficiently small $\kappa>0$. 
Therefore, 
taking a sufficiently small $c>0$ if necessary, 
we apply Theorem~\ref{Theorem:2.4} to see that \eqref{prob} with \eqref{ini} possesses a solution in $[0,\hat{T}_\kappa)$ and
$$
T_\kappa\ge\hat{T}_\kappa=\exp(c\kappa^{-(p-1)})
$$
for all sufficiently small $\kappa>0$. 

In the case of $A=N$, setting 
$$
\check{T}_\kappa:=\exp(c\kappa^{-\frac{p-1}{p}})=\exp(c\kappa^{-\frac{1}{N+1}}), 
$$
similarly to \eqref{eq:4.6} and \eqref{eq:4.7}, we have 
\begin{equation*}
\begin{split}
 & \dashint_{B_+(0,\sigma)}\Phi\left(\check{T}_\kappa^{\frac{N}{2}}\kappa(1+|y|)^{-A}\right)\,dy\\
 & \le C\kappa\check{T}_\kappa^{\frac{N}{2}}(\log\check{T}_\kappa)^N\,\,
\dashint_{B_+(x,\sigma)}\,(1+|y|)^{-N}\,dy
\le C\kappa\check{T}_\kappa^{\frac{N}{2}}\sigma^{-N}(\log\check{T}_\kappa)^{N+1}
=Cc^{N+1}\check{T}_\kappa^{\frac{N}{2}}\sigma^{-N}
\end{split}
\end{equation*}
for $0<\sigma\le\check{T}_\kappa^{1/2}$ and sufficiently small $\kappa>0$. 
Then we apply the same argument as in the case of $A>N$ to see that 
$$
T_\kappa\ge\check{T}_\kappa=\exp(c\kappa^{-\frac{p-1}{p}})
$$
for sufficiently small $\kappa>0$. 
Thus assertion~(i) follows. 

We show assertion~(ii). 
Let $1<p<p_*$ and $0<A<N$. 
Let $c$ be a sufficiently small positive constant and set 
$$
\tilde{T}_\kappa:=c\kappa^{-\left(\frac{1}{2(p-1)}-\frac{A}{2}\right)^{-1}}.
$$
Then 
$$
\sup_{x\in D}\,\dashint_{B_+(x,\tilde{T}_\kappa^{\frac{1}{2}})}\,\kappa\psi(y)\,dy
\le C\kappa\,\dashint_{B_+(0,\tilde{T}_\kappa^{\frac{1}{2}})}\,(1+|y|)^{-A}\,dy
\le C\kappa \tilde{T}_\kappa^{-\frac{A}{2}}
=Cc^{\frac{1}{2(p-1)}-\frac{A}{2}}\tilde{T}_\kappa^{-\frac{1}{2(p-1)}}
$$
for sufficiently small $\kappa>0$. Then we have assertion~(ii) in the case where $1<p<p_*$ and $0<A<N$. 
Similarly, we can prove assertion~(ii) in the other cases and assertion~(ii) follows. 
Thus the proof of Proposition~\ref{Proposition:4.2} is complete.
$\Box$ 
\vspace{5pt}

Finally, we show 
that $\lim_{\kappa\to 0}T_\kappa=\infty$ does not necessarily hold 
for problem~\eqref{prob} if $\psi$ has an exponential growth as $x_N\to\infty$. 
\vspace{5pt}
\\
{\bf Proof of Theorem~\ref{Theorem:1.3}.}  Let $\kappa>0$ and $\delta>0$. 
It follows from Theorem~\ref{Theorem:2.1} that
\begin{equation*}
\begin{split}
\gamma_1T_\kappa^{\frac{N}{2}-\frac{1}{2(p-1)}} 
& > \exp\left(-(1+\delta)\frac{x_N^2}{4T_\kappa}\right)\int_{B_+(x,T_\kappa^{1/2})} \kappa\psi(y)\,dy\\
& \ge C \exp\left(-(1+\delta)\frac{x_N^2}{4T_\kappa}\right) \kappa T_\kappa^{\frac{N}{2}}\exp\biggl(\lambda(x_N-T_\kappa^\frac{1}{2})^2\biggr)\\
& \ge C\kappa T_\kappa^{\frac{N}{2}} \exp\biggl\{ \left( \lambda-\frac{1+\delta}{4T_\kappa}\right)x_N^2\biggr\} 
\exp\biggl( -2\lambda T_\kappa^\frac{1}{2}x_N+\lambda T_\kappa \biggr)
\end{split}
\end{equation*}
for all $x\in D_{T_\kappa}$, where $\gamma_1$ is as in Theorem~\ref{Theorem:2.1}. 
Letting $x_N\to\infty$, we see that $\lambda-(1+\delta)/4T_\kappa\le 0$. 
Since $\delta>0$ is arbitrary, we obtain 
\begin{equation} 
\label{eq:4.8}
\limsup_{\kappa\to +0}\,T_\kappa\le(4\lambda)^{-1}. 
\end{equation}
On the other hand, it follows that 
\begin{equation*}
\begin{split}
\dashint_{B_+(x, \tilde{T}_\delta^{1/2})}
\exp\left(-(1-\delta)\frac{y_N^2}{4\tilde{T}_\delta}\right)\kappa\exp{(\lambda y_N^2)}\,dy=\kappa,
\quad
x\in D_{\tilde{T}_\delta},
\end{split}
\end{equation*}
where $\tilde{T}_\delta := (1-\delta)/4\lambda$. 
Then we deduce from Theorem~\ref{Theorem:2.3} that $T_\kappa\ge \tilde{T}_\delta$ for  sufficiently small $\kappa>0$. 
Since $\delta>0$ is arbitrary, we obtain
$\liminf_{\kappa\to +0}T_\kappa\ge(4\lambda)^{-1}$. 
This together with \eqref{eq:4.8} implies \eqref{eq:1.6}. 
Thus Theorem~\ref{Theorem:1.3} follows.
$\Box$
\vspace{5pt}

\appendix
\def\thesection{}
\section{Appendix}
By Theorem~\ref{Theorem:1.1}, we  obtain Tables~1, 2 and 3. These tables show the behavior of the life span $T(\kappa\psi)$ as $\kappa\to\infty$ when $\psi$ is as in Theorem~\ref{Theorem:1.1}, that is, 
 $$
  \psi(x):=|x|^{-A}\biggr[\log\left(e+\frac{1}{|x|}\right)\biggr]^{-B}\chi_{B_+(0,1)}(x)\in L^1({\bf R}^N_+)
  \setminus L^\infty({\bf R}^N_+), 
  $$
  where $0\le A \le N$ and
\begin{equation*}
B>0\quad\mbox{if}\quad A=0,
\qquad
B\in{\bf R}\quad\mbox{if}\quad 0<A<N,
\qquad
B>1\quad\mbox{if}\quad A=N.
\end{equation*}
For simplicity of notation, we write $T_\kappa$ instead of  $T(\kappa\psi)$.

\begin{table}[htbp]
  \centering
  \caption{The behavior of $T_\kappa$ in the case of $1<p<p_*$ (as $\kappa\to\infty$)}
\renewcommand{\arraystretch}{2.0}
  \begin{tabular}{|c|c|c|} \hline
    \backslashbox{$B$}{$A$} & $A<N$ & $A=N$ \\ \hline
    $B>1$ & $T_\kappa\thicksim\left[\kappa(\log\kappa)^{-B}\right]^{-\frac{2(p-1)}{-A(p-1)+1}}$ & $T_\kappa\thicksim\left[\kappa(\log\kappa)^{-B+1}\right]^{-\frac{2(p-1)}{-A(p-1)+1}}$ \\ \hline
    $B\le1$ & $T_\kappa\thicksim\left[\kappa(\log\kappa)^{-B}\right]^{-\frac{2(p-1)}{-A(p-1)+1}}$ & $0$ \\ \hline
  \end{tabular}
\renewcommand{\arraystretch}{2.0}
\end{table}
\begin{table}[htbp]
  \centering
  \caption{The behavior of $T_\kappa$ in the case of $p>p_*$ (as $\kappa\to\infty$)}
\renewcommand{\arraystretch}{2.0}
  \begin{tabular}{|c|c|c|c|} \hline
    \backslashbox{$B$}{$A$} & $A<\frac{1}{p-1}$ & $A=\frac{1}{p-1}$ & $\frac{1}{p-1} < A \le N$ \\ \hline
    $B>0$ 	& $T_\kappa\thicksim\left[\kappa(\log\kappa)^{-B}\right]^{-\frac{2(p-1)}{-A(p-1)+1}}$ 
    		& $|\log T_\kappa|\thicksim\kappa^\frac{1}{B}$ 
			&  		$0$\\ 
	\hline
    $B=0$ 	& $T_\kappa\thicksim\left[\kappa(\log\kappa)^{-B}\right]^{-\frac{2(p-1)}{-A(p-1)+1}}$ 
    		& $0$ 
    		& $0$ \\
	\hline
    $B<0$ 	& $T_\kappa\thicksim\left[\kappa(\log\kappa)^{-B}\right]^{-\frac{2(p-1)}{-A(p-1)+1}}$ 
    		& $0$ 
			& $0$ \\ 
	\hline
  \end{tabular}
\renewcommand{\arraystretch}{2.0}
\end{table}
\begin{table}[htbp]
   \centering
  \caption{The behavior of $T_\kappa$ in the case of $p=p_*$ (as $\kappa\to\infty$)}
\renewcommand{\arraystretch}{2.0}
  \begin{tabular}{|c|c|c|} \hline
    \backslashbox{$B$}{$A$} & $A<N$ & $A=N$ \\ \hline
    $B>N+1$ & $T_\kappa\thicksim\left[\kappa(\log\kappa)^{-B}\right]^{-\frac{2(p-1)}{-A(p-1)+1}}$ 
    		& $|\log T_\kappa|\thicksim\kappa^\frac{1}{B-N-1}$ \\ 
		\hline
    $B=N+1$ & $T_\kappa\thicksim\left[\kappa(\log\kappa)^{-B}\right]^{-\frac{2(p-1)}{-A(p-1)+1}}$ 
    		& $0$ \\ \hline
    $B<N+1$ & $0$ 
    		& $0$  \\ \hline
  \end{tabular}
\renewcommand{\arraystretch}{2.0}
\end{table}
\newpage
By Theorem~\ref{Theorem:1.2}, we  obtain  Tables~4 and 5. These tables show the behavior of the life span $T(\kappa\psi)$ as $\kappa\to+0$ when $\psi$ is as in Theorem~\ref{Theorem:1.2}, that is,  $\psi(x)=(1+|x|)^{-A} (A>0)$. 

\begin{table}[htbp]
  \centering
  \caption{The behavior of $T_\kappa$ in the case of $A\neq N$ (as $\kappa\to+0$)}
\renewcommand{\arraystretch}{2.0}
  \begin{tabular}{|c|c|c|c|} \hline
    \backslashbox{$p$}{$A$} & $A<\frac{1}{p-1}$ & $A=\frac{1}{p-1}$ & $A>\frac{1}{p-1}$ \\ \hline
    $p<p_*$ 
    & $T_\kappa\thicksim\kappa^{-\left(\frac{1}{2(p-1)}-\frac{1}{2} \min \{A,N\} \right)^{-1}}$ 
    & $T_\kappa\thicksim\kappa^{-\left(\frac{1}{2(p-1)}-\frac{N}{2} \right)^{-1}}$ 
    & $T_\kappa\thicksim\kappa^{-\left(\frac{1}{2(p-1)}-\frac{N}{2} \right)^{-1}}$\\ \hline
    $p=p_*$ 
    & $T_\kappa\thicksim\kappa^{-\left(\frac{1}{2(p-1)}-\frac{A}	{2} \right)^{-1}}$ 
    & ($A=N$, see Table 5)
    & $\log T_\kappa\thicksim \kappa^{-(p-1)}$\\ \hline
    $p>p_*$ 
    & $T_\kappa\thicksim\kappa^{-\left(\frac{1}{2(p-1)}-\frac{A}	{2} \right)^{-1}}$ 
    & $\infty$ 
    & $\infty$\\ \hline
  \end{tabular}
\renewcommand{\arraystretch}{2.0}
 \end{table}
\begin{table}[htbp]
    \centering
  \caption{The behavior of $T_\kappa$ in the case of $A= N$ (as $\kappa\to+0$)}
\renewcommand{\arraystretch}{3.0}
  \begin{tabular}{|c|c|} \hline
    \backslashbox{$p$}{$A$} & $A=N$  \\ \hline
    $p<p_*$ 
    & $T_\kappa\thicksim\displaystyle{\left(\frac{\kappa^{-1}}{\log(\kappa^{-1})}\right)^{\left(\frac{1}{2(p-1)}-\frac{N}{2}\right)^{-1}}}$ 
    \\ \hline
    $p=p_*$ 
    & $\log T_\kappa\thicksim \kappa^{-\frac{p-1}{p}}$ 
    \\ \hline
    $p>p_*$ 
    & $\infty$
    \\ \hline
  \end{tabular}
\renewcommand{\arraystretch}{2.0}
  \end{table}
\vspace{5pt}
 
\newpage

\vspace{5pt}
\noindent
{\bf Acknowledgments.} 
The author of this paper is grateful to Professor K. Ishige for mathematical discussions and proofreading the manuscript. 
The author of this paper is grateful to Professor S. Okabe for carefully proofreading of the manuscript. 
Finally,
the author of this paper would like to thank the referees for carefully reading the manuscript and relevant remarks. 



\begin{thebibliography}{10}

\bibitem{ACR1} 
J.~M. Arrieta, A.~N. Carvalho and A. Rodr\'{\i}guez-Bernal, 
Parabolic problems with nonlinear boundary conditions and critical nonlinearities, 
J. Differ. Equ. {\bf 156} (1999), 376--406.

\bibitem{AR}
J.~M.~Arrieta and A.~Rodr\'{\i}guez-Bernal, 
Non well posedness of parabolic equations with supercritical nonlinearities, 
Commun. Contemp. Math. {\bf 6} (2004), 733--764.


\bibitem{DFL} 
K.~Deng, M.~Fila and H.~A.~Levine, 
On critical exponents for a system of heat equations coupled in the boundary conditions, 
Acta Math. Univ. Comenianae {\bf 63} (1994), 169--192. 

\bibitem{D}
F.~Dickstein, Blowup stability of solutions of the nonlinear heat equation with a large life span. J. Differ. Equ. {\bf 223}, 303--328 (2006)


\bibitem{FR}
J.~Fern\'andez Bonder and J.~D.~Rossi, 
Life span for solutions of the heat equation with a nonlinear boundary condition, Tsukuba J. Math. {\bf 25} (2001), 215--220.

\bibitem{FiloK}
J.~Filo and J.~Ka\v{c}ur,
Local existence of general nonlinear parabolic systems.
Nonlinear Anal. {\bf 24} (1995), 1597--1618.

\bibitem{GL} 
V.~A.~Galaktionov and H.~A.~Levine, 
On critical Fujita exponents for heat equations with nonlinear flux conditions on the boundary, 
Israel J. Math. {\bf 94} (1996), 125--146.

\bibitem{GW}
C.~Gui and X.~Wang, Life spans of solutions of the Cauchy problem for a semilinear heat equation, J. Differ. Equ. {\bf 115} (1995), 166--172

\bibitem{HI18} 
K.~Hisa and K.~Ishige,
Existence of solutions for a fractional semilinear parabolic equation 
with singular initial data, 
Nonlinear Anal. {\bf 175} (2018), 108--132. 

\bibitem{HI19}
K.~Hisa and K.~Ishige,
Solvability of the heat equation with a nonlinear boundary condition, 
SIAM J. Math. Anal. {\bf 51} (2019), 565--594.


\bibitem{IkSo1}
M.~Ikeda and M.~Sobajima, 
Sharp upper bound for lifespan of solutions to some critical semilinear parabolic, dispersive and hyperbolic equations via a test function method, Nonlinear Anal. {\bf 182} (2019), 57--74.

\bibitem{IK}
K.~Ishige and T.~Kawakami,
Global solutions of the heat equation with a nonlinear boundary condition,
Calc. Var. Partial Differential Equations {\bf 39} (2010), 429--457.

\bibitem{IS1}
K.~Ishige and R.~Sato, 
Heat equation with a nonlinear boundary condition and uniformly local $L^r$ spaces,
Discrete Contin. Dyn. Syst. {\bf 36} (2016), 2627--2652.

\bibitem{IS2}
K.~Ishige and R.~Sato, 
Heat equation with a nonlinear boundary condition and growing initial data,
Differential Integral Equations {\bf 30} (2017), 481--504. 

\bibitem{LN}
T.~Y.~Lee and W.~M.~Ni, 
Global existence, large time behavior and life span of
solutions of a semilinear parabolic Cauchy problem, 
Trans. Amer. Math. Soc. {\bf 333} (1992), 365--378.

\bibitem{MY1}
N.~Mizoguchi and E.~Yanagida, 
Blowup and life span of solutions for a semilinear parabolic equation, 
SIAM J. Math. Anal. {\bf 29} (1998),  1434--1446. 

\bibitem{MY2}
N.~Mizoguchi and E.~Yanagida, 
Life span of solutions with large initial data in a semilinear parabolic equation, 
Indiana Univ. Math. J. {\bf 50} (2001),  591--610.

\bibitem{OY}
T.~Ozawa and Y.~Yamauchi, Life span of positive solutions for a semilinear heat equation with general non-decaying initial data, J. Math. Anal. Appl. {\bf 379} (2011), 518--523.

\bibitem{S}
S.~Sato, Life span of solutions with large initial data for a superlinear heat equation, J. Math. Anal. Appl. {\bf 343} (2008), 1061--1074.

\bibitem{YY}
M.~Yamaguchi and Y.~Yamauchi, Life span of positive solutions for a semilinear heat equation with non-decaying initial data, Differential Integral Equations {\bf 23} (2010), 1151--1157. 

\bibitem{Y}
Y.~Yamauchi, Life span of solutions for a semilinear heat equation with initial data having positive limit inferior at infinity, Nonlinear Anal. {\bf 74} (2011),  5008--5014.


\end{thebibliography}
\end{document}